\documentclass[10pt,a4paper]{article}
\usepackage{cite}

\usepackage{CJK,fancybox,fancyhdr,color,graphicx,amsmath,amsthm,amssymb,enumerate,amsfonts}
\usepackage{mathrsfs}
\usepackage{indentfirst,latexsym,bm}
\usepackage{enumerate,mathdots}
\usepackage{blkarray}
\usepackage{mathtools}
\usepackage{float}
\usepackage{caption}
\usepackage{subfigure}
\usepackage{anysize}
\usepackage{marvosym,titling}
\usepackage{authblk}
\usepackage[hang]{footmisc}
\usepackage[T1]{fontenc}
\usepackage{algorithmicx,algorithm, algpseudocode}
\usepackage{xcolor}%
\usepackage{colortbl,booktabs}

\marginsize{2cm}{2cm}{2.54cm}{2cm} 

\newtheorem{thm}{Theorem}[section]
\newtheorem{cor}[thm]{Corollary}
\newtheorem{lem}[thm]{Lemma}

\newtheorem{defn}[thm]{Definition}

\numberwithin{equation}{section}

\renewcommand{\proof}{\noindent\textbf{Proof.}}

\title{The standard form and convergence theory of the relaxation Kaczmarz-Tanabe method for solving linear systems}

\date{}
\makeatletter
\newenvironment{breakablealgorithm}
{
\begin{center}
\refstepcounter{algorithm}
\hrule height.8pt depth0pt \kern2pt
\renewcommand{\caption}[2][\relax]{
{\raggedright\textbf{\ALG@name~\thealgorithm} ##2\par}%
\ifx\relax##1\relax 
\addcontentsline{loa}{algorithm}{\protect\numberline{\thealgorithm}##2}%
\else 
\addcontentsline{loa}{algorithm}{\protect\numberline{\thealgorithm}##1}%
\fi
\kern2pt\hrule\kern2pt
}
}{
\kern2pt\hrule\relax
\end{center}
}

\makeatother

\begin{document}


%

\author{Chuan-gang Kang\\School of Mathematical Sciences, Tiangong University, Tianjin, 300387, People's Republic of China\vspace{-0.5cm}}
\footnotetext{Email address: ckangtj@tiangong.edu.cn}
\maketitle

\begin{abstract} The Kaczmarz method is a popular iterative method for solving consistent, overdetermined linear system such as  medical imaging in computerized tomography. The Kaczmarz's iteration repeatedly scans all equations in order, which leads to lower computational efficiency especially in solving a large scale problem. The standard form of Kaczmarz-Tanabe's iteration proposed recently effectively overcomes the computational redundancy problem of the Kaczmarz method. In this paper, we introduce relaxation parameters ${\bf u}=(\mu_1,\ldots,\mu_m)$ into the Kaczmarz-Tanabe method based on the relaxation Kaczmarz method, and consider the standard form and convergence of this combination. Moreover, we analyze and prove the sufficient conditions for convergence of the relaxation Kaczmarz-Tanabe method, i.e., $\mu_i\in (0,2)$. Numerical experiments show the convergence characteristics of the relaxation Kaczmarz-Tanabe method corresponding to these parameters.
\end{abstract}




\section{Introduction}
\label{section.introduction}

In medical imaging tomography (see, i.e., \cite{Ledley1988,Herman2010Fundamentals,Natterer1986}) such as computerized tomography (CT), magnetic resonance imaging (MRI) and electrical impedance tomography (EIT), people are often required to solve the following linear system of equations \cite{Jinbt2020Adaptive,Alberti2016The,Chow2014},
\begin{align}\label{linear.system}
  Ax=b,
\end{align}
where $A\in\mathbb{R}^{m\times n}, b\in\mathbb{R}^m$ are also called projection matrix and measurement vector, respectively. We suppose that \eqref{linear.system} is over-determined and consistent, and $x^*$ is a true solution. When \eqref{linear.system} has many solutions, $x^\dagger$ is used to denote the minimum least-squares solution.

The Kaczmarz method \cite{Kaczmarz1937} is one of the most popular methods to solve \eqref{linear.system} in computerized tomography. Let $A=(a_1,a_2,\ldots,a_m)^T$, Kaczmarz's iteration is written as
\begin{align}\label{Kaczmarz.iteration}
  x_k=x_{k-1}+\frac{b_i-\langle a_i,x_{k-1}\rangle}{\|a_i\|_2^2}a_i,\quad k=1,2,\ldots,
\end{align}
where $i=\bmod(k-1,m)+1$, $\langle x,y\rangle=x^Ty$ and $\|x\|_2=\sqrt{\langle x,x\rangle}$ denote the inner product of $x,y$ and the square norm of $x$ in $\mathbb{R}^n$, respectively.

The Kaczmarz method was originally discovered by Kaczmarz in 1937, and it was rediscovered by Gordon, Bender and Herman \cite{Gordon1970} to deal with X-ray photography as an algebraic reconstruction technique (ART) in 1970. The convergence of the Kaczmarz method was established by Tanabe in 1971. Recently, Kang and Zhou gave the convergence rate results of the Kaczmarz method \cite{Kang2021The}. In order to improve computing speed, Eggermont \cite{Eggermont1981}, Censor \cite{Censor1983strong}, et al. considered the problem of block calculation and adopted approximate and executable iterative forms of the block Kaczmarz method to solve this problem. The strictly block Kaczmarz method was given by Rebrova and Needell \cite{Rebrova2021}. In addition, the simultaneous iterative reconstructive technique (SIRT) as another important branch is also developed rapidly, and there appeared many known methods) such as Cimmino method\cite{Cimmino1938}, Landweber method \cite{Landweber51}, CAV method \cite{Censor2001Component}, DROP method \cite{Censor2007On}, and SART method\cite{Hansen2018AIRtools}, etc.

With the development of hardware technology and computing technology, the idea of stochastic theory is continuously integrated into the non-randomized methods. In this context, people introduced and investigated the randomized Kaczmarz method \cite{Strohmer2009,Jiao2017preasymptotic,Needell2010,Needell2014paved,Needell2015Randomized} and greedy randomized Kaczmarz methods \cite{Bai2018onconvergence,Bai2018ongreedy}.

From the perspective of SIRT methods, Kang \cite{Kang2021Convergence} considered the matrix-vector form of the Kaczmarz method, i.e., the Kaczmarz-Tanabe's iteration which is named by Popa \cite{Popa2018}, and the iteration reads
\begin{align}\label{Kaczmarz.Tanabe.iteration.orig}
  y_{k+1}=(I-A_{\mathcal{S}}^TMA)y_k+A_{\mathcal{S}}^TMb,\qquad k=0,1,2,\ldots,
\end{align}
where $I$ denotes identity matrix of whatever size appropriate to the context, and
\begin{align*}
  &P_i=I-\frac{a_ia_i^T}{\|a_i\|_2^2},\quad
  Q_m=I, Q_j=P_mP_{m-1}\ldots P_{j+1}, \quad
  Q=P_mP_{m-1}\cdots P_1,\\
  &M=\text{diag}(1/\|a_1\|_2^2,1/\|a_2\|_2^2,\ldots,1/\|a_m\|_2^2),\quad
  A_\mathcal{S}=(Q_1a_1,Q_2a_2,\ldots,Q_ma_m)^T.
\end{align*}

Compared with the classic form of the Kaczmarz-Tanabe iteration, the iteration \eqref{Kaczmarz.Tanabe.iteration.orig} has been greatly improved. However, each column of $A_{\mathcal{S}}$ is the product of matrix $Q_i$ and column vector $a_i$ from $A$. Notice this defect of $A_{\mathcal{S}}$, Kang \cite{Kang2023Thestandard} considered the decomposition of $A_{\mathcal{S}}$ and gave the standard form of the Kaczmarz-Tanabe's iteration, i.e.,
\begin{align}\label{Kaczmarz.Tanabe.iteration.standard}
  y_{k+1}=y_k+A^TC^TM(b-Ay_k), \qquad k=0,1,2,\ldots,
\end{align}
where $C$ is defined as \cite[Corollary 2.7]{Kang2023Thestandard}. The iterative scheme \eqref{Kaczmarz.Tanabe.iteration.standard} has no essential difference from the iterative formula of the SIRT methods.

In this paper, we combine the relaxation Kaczmarz method and the Kaczmarz-Tanabe method, and introduce relaxation parameters into Kaczmarz-Tanabe method and propose the relaxation Kaczmarz-Tanabe method. Just as the convergence of the relaxation method is related to the relaxation parameters, the convergence of the relaxation Kaczmarz-Tanabe is also related to the relaxation parameters. Hence we further consider the convergence conditions (about the relaxation parameters) and convergence results of the relaxation Kaczmarz-Tanabe method.

The rest of this work is organized as follows. In Section \ref{section.relaxation.Kaczmarz.method}, we consider the relaxation Kaczmarz method and derive the sufficient conditions for its convergence. In Section \ref{section.relaxation.Kaczmarz-Tanabe.method}, we consider the standard form of the relaxation Kaczmarz-Tanabe method, and discuss their properties such as decomposition of $A_{\mathcal{S}}({\bf u})$, convergence condition, convergence, and so on. Section \ref{section.algorithm} is about the algorithm flow to compute $C({\bf u})$ appearing in \eqref{definition.Omega.u} and \eqref{Cu.decomposition.formula}. In Section \ref{section.numerical.tests}, we verify the effectiveness of the Kaczmarz-Tanabe method with given relaxation parameters by several classical numerical examples.

\section{The relaxation Kaczmarz method and its convergence}
\label{section.relaxation.Kaczmarz.method}

The relaxation Kaczmarz's iteration is written as
\begin{align}\label{relaxation.kaczmarz.iteration}
  x_k=x_{k-1}+\mu_i\frac{b_i-\langle a_i, x_{k-1}\rangle}{\|a_i\|_2^2}a_i, \quad k=1,2,\ldots,
\end{align}
where $i=\bmod(k-1,m)+1$ and $\mu_i$ are relaxation parameters. Let
\begin{align}\label{Pi.mu}
  P_i(\mu_i)=I-\mu_i\frac{a_ia_i^T}{\|a_i\|_2^2}, \quad i=1,2,\ldots,m,
\end{align}
then it follows from \eqref{relaxation.kaczmarz.iteration} that
\begin{align}\label{relaxation.kaczmarz.iteration.2}
  x_k=P_i(\mu_i)x_{k-1}+\mu_i\frac{b_i}{\|a_i\|_2^2}a_i,\quad k=1,2,\ldots.
\end{align}
Denote
  \begin{align}\label{Pi.mu.3}
    P_i(\mu):=I-\mu a_ia_i^T/\|a_i\|_2^2,
  \end{align}
then the following lemma gives a sufficient condition for the nonexpansion of linear operator $P_i(\mu)$.
\begin{lem}\label{lemma.Pi.norm} Let $\{P_i(\mu), i=1,2,\cdots,m\}$ be defined as in \eqref{Pi.mu.3}, then, for any $0<\mu<2$,
\begin{align*}
  \|P_i(\mu)x\|_2\le\|x\|_2,
\end{align*}
hold and the equality hold iif $x\in N(a_i^T)\triangleq\{x: a_i^Tx=0\}$. Moreover, for any $\mu\in \mathbb{R}$ there also holds $P_i(\mu)x=x$ when $x\in N(a_i^T)$.
\end{lem}
\proof \quad For any $1\le i\le m$, we have
\begin{align}\label{Pi.mu.2}
  \|P_i(\mu)x\|_2^2=\|(I-\mu\frac{a_ia_i^T}{\|a_i\|_2^2})x\|_2^2=\|x\|_2^2-\mu(2-\mu)\|a_i^Tx\|_2^2/\|a_i\|_2^2.
\end{align}
If $x\notin N(a_i^T)$ and $0<\mu<2$, then $\mu(2-\mu)\|a_i^Tx\|_2^2>0$, consequently, $\|P_i(\mu)x\|_2<\|x\|_2$. Conversely, if $x\in N(a_i^T)$, then $a_i^Tx=0$, and it follows from \eqref{Pi.mu.2} that $\|P_i(\mu)x\|_2^2=\|x\|_2^2$. Moreover, there also holds $P_i(\mu)x=x$ when $x\in N(a_i^T)$. \qed

Let $\{x_k,k>0\}$ be generated by \eqref{relaxation.kaczmarz.iteration} and $e_k=x_k-P_{N(A)}x_0-x^\dagger$, thus
\begin{align}\label{relaxation.kaczmarz.error}
    e_k=(I-\mu_i\frac{a_ia_i^T}{\|a_i\|_2^2})e_{k-1}=P_i(\mu_i)e_{k-1},
\end{align}
where $i=\bmod(k-1,m)+1$, and we have the following set property.
\begin{lem}\label{lem.set.property.relaxation.Kaczmarz} Let $\{x_k,k>0\}$ be generated by \eqref{relaxation.kaczmarz.iteration}, then, for any $\mu_i\in \mathbb{R}$,
\begin{align*}
  e_k\in N(A)^{\bot}
\end{align*}
holds.
\end{lem}
\proof \quad We prove the result by mathematical induction. First, we know that $e_0=y_0-x^\dagger-P_{N(A)}y_0\in N(A)^{\bot}$. Second, we suppose that for any given $k\ge 0$ there holds $e_k\in N(A)^{\bot}$, then $e_{k+1}\in N(A)^{\bot}$. In fact, for any $z\in N(A)$,
\begin{align*}
  \langle e_{k+1},z\rangle=\langle P_i(\mu_i)e_k,z\rangle=\langle e_k,z\rangle-\mu_i\frac{b_i}{\|a_i\|_2^2}\langle a_i,z\rangle=0.
\end{align*}
This proves $e_{k+1}\in N(A)^{\bot}$. To sum up, $e_k\in N(A)^{\bot}$ for any $k\ge 0$.\qed

Lemma \ref{lem.set.property.relaxation.Kaczmarz} shows that, for any relaxation parameter $\mu_i$, the iterative error $\{e_k,k\ge 0\}$ of the relaxation Kaczmarz method belong to $N(A)^{\bot}$. The following theorem gives the convergence result of the relaxation Kaczmarz method for solving a consistent linear system when $0<\mu_i<2$ (where $i=1,2,\ldots,m$).
\begin{thm}\label{thm.relaxation.kaczmarz.convergence} For a consistent linear system $Ax=b$, let $\{x_k,k\ge 0\}$ be generated by \eqref{relaxation.kaczmarz.iteration} and relaxation parameter $0<\mu_i<2$. Then the sequence $\{x_k,k\ge 0\}$ converges to the solution $x^\dagger$ when $x_0\in R(A^T)$.
\end{thm}
\proof \quad Let $\bar{P}_i(\mu_i)\triangleq P_i(\mu_i)|_{N(A)^{\bot}}$, when $0<\mu_i<2$, then from Lemma \ref{lemma.Pi.norm} there holds $\|\bar{P}_i(\mu_i)\|_2<1$, consequently $\zeta\triangleq\max\limits_{1\le i\le m}\|\bar{P}_i(\mu_i)\|<1$. It follows from \eqref{relaxation.kaczmarz.error} that
\begin{align*}
  \|e_k\|_2=\|P_i(\mu_i)e_{k-1}\|_2.
\end{align*}
By Lemma \ref{lem.set.property.relaxation.Kaczmarz}, then
\begin{align*}
  \|e_k\|_2=\|\bar{P}_i(\mu_i)e_{k-1}\|_2\le \|\bar{P}_i(\mu_i)\|\|e_{k-1}\|_2<\zeta\|e_{k-1}\|_2<\ldots<\zeta^k\|e_0\|_2,
\end{align*}
which means $\|e_k\|_2\rightarrow 0$ as $k\rightarrow \infty$, i.e. $x_k$ converges to $x^\dagger+P_{N(A)}x_0$. Moreover $x_k$ also converges to $x^\dagger$ when $x_0\in R(A^T)$. \qed

\section{The standard form of the relaxation Kaczmarz-Tanabe method}
\label{section.relaxation.Kaczmarz-Tanabe.method}
We start with the first equation of \eqref{linear.system} and perform the iteration \eqref{relaxation.kaczmarz.iteration.2} in turn until the $m$-th equation,
and denote the $m$-th iteration with $x_m$, then get
\begin{align}\label{rexlaxation.Kaczmarz.Tanabe.iteration}
  x_m&=P_m(\mu_m)\cdots P_1(\mu_1)x_0\nonumber\\
     &\quad +\mu_1P_m(\mu_m)\cdots P_2(\mu_2)\frac{b_1}{\|a_1\|_2^2}a_1+\mu_2P_m(\mu_m)\cdots P_3(\mu_3)\frac{b_2}{\|a_2\|_2^2}a_2+\cdots+\mu_m\frac{b_m}{\|a_m\|_2^2}a_m.
\end{align}
For $1\le j\le m-1$, let
\begin{align}
  &{\bf u}:=(\mu_1,\ldots,\mu_m),\quad {\bf u}_j:=(\mu_{j+1},\ldots,\mu_m),\label{definition.u.uj}\\
  &Q({\bf u}):=P_m(\mu_m)\cdots P_1(\mu_1), \quad Q_j({\bf u}_j):=P_m(\mu_m)\cdots P_{j+1}(\mu_{j+1}).\label{definition.Qmu}
\end{align}
If $\mu_i\equiv 1$ for each $1\le i\le m$, then $Q({\bf u})=Q, Q({\bf u}_j)=Q_j$. In addition, we agree that $Q_m({\bf u}_m)=I$, where $I$ denote identity matrix of whatever size appropriate to the context.

Let
\begin{align}\label{definition.AS.mu.and.Lambda}
  A_{\mathcal{S}}({\bf u})\triangleq(Q_1({\bf u}_1)a_1,Q_2({\bf u}_2)a_2,\ldots,Q_m({\bf u}_m)a_m\big)^T, \quad \Lambda=\mathbf{diag}(\mu_1,\ldots,\mu_m).
\end{align}
Thus we get from \eqref{rexlaxation.Kaczmarz.Tanabe.iteration}, \eqref{definition.u.uj}, \eqref{definition.Qmu} and \eqref{definition.AS.mu.and.Lambda} that
\begin{align*}
  x_m=Q({\bf u})x_0+[A_{\mathcal{S}}({\bf u})]^T\Lambda Mb.
\end{align*}
Consequently, it follows by repeating this iteration that
\begin{align}\label{relaxation.Kaczmarz.Tanabe.iteration.3}
  x_{k\cdot m}=Q({\bf u})x_{(k-1)\cdot m}+[A_{\mathcal{S}}({\bf u})]^T\Lambda Mb, \quad k=1,2,\ldots.
\end{align}
Taking $y_k=x_{k\cdot m}$, then we get the \textbf{relaxation Kaczmarz-Tanabe's iteration}
\begin{align}\label{relaxation.Kaczmarz.Tanabe.iteration.AS}
  y_k=Q({\bf u})y_{k-1}+[A_{\mathcal{S}}({\bf u})]^T\Lambda Mb, \quad k=1,2,\ldots,
\end{align}
where $y_0=x_0$.

\begin{lem} \label{lem.Qmu}Let $Q(\bf{u})$ and $Q_j({\bf u}_j)$ be defined as in \eqref{definition.Qmu}, then there holds
\begin{align}\label{property.Qmu}
  Q({\bf u})=I-[A_{\mathcal{S}}({\bf u})]^T\Lambda MA.
\end{align}
\end{lem}
\proof \quad First, we have from \eqref{definition.Qmu} that
\begin{align*}
  Q({\bf u})&=Q_1({\bf u}_1)-\mu_1 Q_1({\bf u}_1)\frac{a_1a_1^T}{\|a_1\|_2^2}.
\end{align*}
By $Q_{j-1}({\bf u}_{j-1})=Q_j({\bf u}_j)-\mu_jQ_j({\bf u}_j)\frac{a_ja_j^T}{\|a_j\|_2^2}$, then
\begin{align*}
  Q({\bf u})&=Q_1({\bf u}_1)-\mu_1 Q_1({\bf u}_1)\frac{a_1a_1^T}{\|a_1\|_2^2}\\
        &=Q_2({\bf u}_2)-\big(\mu_1Q_1({\bf u}_1)\frac{a_1a_1^T}{\|a_1\|_2^2}+\mu_2Q_2({\bf u}_2)\frac{a_2a_2^T}{\|a_2\|_2^2}\big)\\
        &=\cdots\\
        &=Q_m({\bf u}_m)- \big(\mu_1Q_1({\bf u}_1)\frac{a_1a_1^T}{\|a_1\|_2^2}+\mu_2Q_2({\bf u}_2)\frac{a_2a_2^T}{\|a_2\|_2^2}+\cdots+\mu_mQ_m({\bf u}_m)\frac{a_ma_m^T}{\|a_m\|_2^2}\big)\\
        &=I-\big(\mu_1Q_1({\bf u}_1)a_1,\mu_2Q_2({\bf u}_2)a_2,\cdots,\mu_mQ_m({\bf u}_m)a_m\big)\mathbf{diag}(\frac{1}{\|a_1\|_2^2},\cdots,\frac{1}{\|a_m\|_2^2})\big(a_1,a_2,\cdots,a_m\big)^T.
\end{align*}
holds. This proves \eqref{property.Qmu}.\qed

By Lemma \ref{lem.Qmu}, we can get the alternative form of \eqref{relaxation.Kaczmarz.Tanabe.iteration.AS}, i.e.,
\begin{align*}
  y_k=y_{k-1}+[A_{\mathcal{S}}({\bf u})]^T\Lambda M(b-Ay_{k-1}), \quad k=1,2,\ldots.
\end{align*}
For the further discussion to \eqref{relaxation.Kaczmarz.Tanabe.iteration.AS}, we extend the conceptions of sequential projection matrix, sequential projection matrix set and sequential compatible to $Q_j({\bf u}_j)$ (see \cite[Definition 2.1 \& Definition 2.2]{Kang2023Thestandard}).

\begin{defn}\label{definition.S_sp} We call $Q_j({\bf u}_j)$ a \textbf{sequential projection matrix with parameters} on $a_{j+1},\ldots,a_m$, and denote the sequential projection matrix set with $S_{sp}(a_1,\ldots,a_m;{\bf u})$, i.e.,
\begin{align}\label{sequential.projection.matrix.set}
  S_{sp}(a_1,a_2,\ldots,a_m;{\bf u})=\{Q_1({\bf u}_1),Q_2({\bf u}_2),\ldots,Q_{m-1}({\bf u}_{m-1})\}.
\end{align}
\end{defn}

\begin{defn}\label{definition.compatible} For any $Q_i({\bf u}_i)\in S_{sp}(a_1,a_2,\ldots,a_m;{\bf u})$, if there exist $\zeta_{i,1},\ldots, \zeta_{i,m}$ such that
\begin{align}\label{sequential.compatible}
  Q_i({\bf u}_i)a_i=\zeta_{i,1}({\bf u}_i)a_1+\ldots \zeta_{i,m}({\bf u}_i)a_m,
\end{align}
then we call $A$ and $S_{sp}(a_1,a_2,\ldots,a_m;{\bf u})$ \textbf{sequentially compatible}. In general, for any $1\le i\le m, i\le j\le m$, if there exist $\zeta_1^{(i,j)}({\bf u}_j),\ldots,\zeta_m^{(i,j)}({\bf u}_j)$ such that
\begin{align}\label{forward.sequential.compatible}
  Q_j({\bf u}_j)a_i=\zeta_1^{(i,j)}({\bf u}_j)a_1+\ldots \zeta_m^{(i,j)}({\bf u}_j)a_m,
\end{align}
then we call $A$ and $S_{sp}(a_1,a_2,\ldots,a_m;{\bf u})$ \textbf{forward sequential compatible}, and call $(\zeta_1^{(i,j)}({\bf u}_j),\ldots,\zeta_m^{(i,j)}({\bf u}_j))$ \textbf{compatible vector} of $Q_j({\bf u}_j)a_i$ on $A$.
\end{defn}

\begin{thm}\label{thm.forward.sequential.compatible} Suppose $A$ has no zero row, and $S_{sp}(a_1,\ldots,a_m;{\bf u})$ is defined by \eqref{sequential.projection.matrix.set}, then $A$ and $S_{sp}(a_1,\ldots,a_m;{\bf u})$ are forward sequential compatible.
\end{thm}
\proof \quad Similar to \cite[Theorem 2.6]{Kang2023Thestandard}, we will complete the proof by mathematical induction. We take the subscript $(i,j)$ of $Q_j({\bf u}_j)a_i$ as an ordered array. Obviously, $a_m^T[Q_m({\bf u}_m)]^T=a_m^T$, that is, \eqref{forward.sequential.compatible} holds for $(i,j)=(m,m)$ and $(\zeta_1^{(m,m)}({\bf u}_m),\ldots,\zeta_m^{(m,m)}({\bf u}_m))=(0,\ldots,0,1)$. In fact, for any $1\le i\le m$, \eqref{forward.sequential.compatible} obviously holds for $a_i^T[Q_m({\bf u}_m)]^T$ because $Q_m=I$ is the agreement. Consequently, As the first step of mathematical induction, we prove \eqref{forward.sequential.compatible} holds for $(i,j)=(m-1,m-1)$.
Actually,
\begin{align*}
  a_{m-1}^T[Q_{m-1}({\bf u}_{m-1})]^T=a_{m-1}^T[P_m(\mu_m)]^T[Q_m({\bf u}_m)]^T=a_{m-1}^T-\mu_m\frac{a_{m-1}^T a_m}{\|a_m\|_2^2}a_m^T.
\end{align*}
So \eqref{sequential.compatible} holds for $i=m-1$, where
\[(\zeta_1^{(m-1,m-1)}({\bf u}_{m-1}),\ldots,\zeta_m^{(m-1,m-1)}({\bf u}_{m-1}))=(0,\ldots,0,1,-\mu_ma_{m-1}^Ta_m/\|a_m\|_2^2).\]

Second, we suppose \eqref{forward.sequential.compatible} holds for any given $(i,j)$ subjected to $s\le i<m$ and $s\le t<j<m$, i.e.,
there exists $(\zeta_1^{(i,j)}({\bf u}_j),\ldots,\zeta_m^{(i,j)}({\bf u}_j))$ such that
\begin{align*}
  a_i^T[Q_j({\bf u}_j)]^T=\zeta_1^{(i,j)}({\bf u}_j)a_1^T+\ldots \zeta_m^{(i,j)}({\bf u}_j)a_m^T.
\end{align*}
Then, we prove \eqref{forward.sequential.compatible} holds for $(i,j)=(s,t)$. Because of $Q_t({\bf u}_t)=Q_{t+1}({\bf u}_{t+1})P_{t+1}(\mu_{t+1})$, then
\begin{align}\label{math.induction.third}
  a_s^T[Q_t({\bf u}_t)]^T=a_s^T[P_{t+1}(\mu_{t+1})]^T[Q_{t+1}({\bf u}_{t+1})]^T=a_s^T[Q_{t+1}({\bf u}_{t+1})]^T-\mu_{t+1}\frac{a_s^Ta_{t+1}}{\|a_{t+1}\|_2^2}a_{t+1}^T[Q_{t+1}({\bf u}_{t+1})]^T.
\end{align}
From the hypothesis, there exist
$(\zeta_1^{(s,t+1)}({\bf u}_{t+1}),\ldots,\zeta_m^{(s,t+1)}({\bf u}_{t+1}))$ and $(\zeta_1^{(t+1,t+1)}({\bf u}_{t+1}),\ldots,\zeta_m^{(t+1,t+1)}({\bf u}_{t+1}))$ such that
\begin{align*}
  &a_s^T[Q_{t+1}({\bf u}_{t+1})]^T=\zeta_1^{(s,t+1)}({\bf u}_{t+1})a_1^T+\ldots+\zeta_m^{(s,t+1)}({\bf u}_{t+1})a_m^T,\\
  &a_{t+1}^T[Q_{t+1}({\bf u}_{t+1})]^T=\zeta_1^{(t+1,t+1)}({\bf u}_{t+1})a_1^T+\ldots+\zeta_m^{(t+1,t+1)}({\bf u}_{t+1})a_m^T.
\end{align*}
Then, it follows from \eqref{math.induction.third} that
\begin{align*}
  &a_s^T[Q_t({\bf u}_t)]^T\\
  &=\zeta_1^{(s,t+1)}({\bf u}_{t+1})a_1^T+\ldots+\zeta_m^{(s,t+1)}({\bf u}_{t+1})a_m^T-\frac{a_s^Ta_{t+1}}{\|a_{t+1}\|_2^2}\big(\zeta_1^{(t+1,t+1)}({\bf u}_{t+1})a_1^T+\ldots+\zeta_m^{(t+1,t+1)}({\bf u}_{t+1})a_m^T\big)\\
  &=\big(\zeta_1^{(s,t+1)}({\bf u}_{t+1})-\frac{a_s^Ta_{t+1}}{\|a_{t+1}\|_2^2}\zeta_1^{(t+1,t+1)}({\bf u}_{t+1})\big)a_1^T+\ldots+\big(\zeta_m^{(s,t+1)}({\bf u}_{t+1})-\frac{a_s^Ta_{t+1}}{\|a_{t+1}\|_2^2}\zeta_m^{(t+1,t+1)}({\bf u}_{t+1})\big)a_m^T.
\end{align*}
Denote
\begin{align*}
  &\big(\zeta_1^{(s,t)}({\bf u}_t),\ldots,\zeta_m^{(s,t)}({\bf u}_t)\big)\\
  &=\big(\zeta_1^{(s,t+1)}({\bf u}_{t+1})-\frac{a_s^Ta_{t+1}}{\|a_{t+1}\|_2^2}\zeta_1^{(t+1,t+1)}({\bf u}_{t+1}),\ldots,\zeta_m^{(s,t+1)}({\bf u}_{t+1})-\frac{a_s^Ta_{t+1}}{\|a_{t+1}\|_2^2}\zeta_m^{(t+1,t+1)}({\bf u}_{t+1})\big).
\end{align*}
This proves \eqref{forward.sequential.compatible} holds for $(i,j)=(s,t)$.

To sum up the above, the conclusion is proved for all $(i,j)$ with respect to $1\le i\le m,i\le j\le m$. Namely, $A$ and $S_{sp}(a_1,a_2,\ldots,a_m;{\bf u})$ generated by relaxation Kaczmarz's iteration are forward sequential compatible.  \qed

In order to obtain the standard form (i.e., matrix-vector form) of the relaxation Kaczmarz-Tanabe's iteration, we first import the following decomposition theorem about $A_{\mathcal{S}}({\bf u})$.

\begin{thm}\label{thm.AS.decomposition.relaxation} Suppose $A$ has no zero rows, and let $A_{\mathcal{S}}({\bf u})$ be defined by \eqref{definition.AS.mu.and.Lambda}, then there exists a unit upper triangular matrix $C({\bf u})\in \mathbb{R}^{m\times m}$ such that $A_\mathcal{S}({\bf u})=C({\bf u})A$.
\end{thm}
\proof \quad According to $A_\mathcal{S}({\bf u})=(Q_1({\bf u}_1) a_1,\ldots,Q_m({\bf u}_m) a_m)^T$ and Theorem \ref{thm.forward.sequential.compatible}, the corollary can be proved by taking $[C({\bf u})](i,j)=\zeta_{i,j}({\bf u}_i)$, where $\zeta_{i,j}({\bf u}_i)$ is determined by the forward sequential compatible of $A$ and $S_{sp}(a_1,\ldots,a_m;{\bf u})$. \qed

Theorem \ref{thm.AS.decomposition.relaxation} is a general result, but it shows that there exists a structure-reserved decomposition in relaxation Kaczmarz-Tanabe's iteration when a relaxation parameter is introduced into Kaczmarz's iteration. However, the decomposition is not unique when $A$ is not row full rank.

By Theorem \ref{thm.AS.decomposition.relaxation}, then we get \textbf{the standard form of the relaxation Kaczmarz-Tanabe's iteration} from \eqref{relaxation.Kaczmarz.Tanabe.iteration.AS} that
\begin{align}\label{normal.form.relaxation.Kaczmarz.Tanabe.iteration}
  y_k=y_{k-1}+A^T[C({\bf u})]^T\Lambda M(b-Ay_{k-1}), \quad k=1,2,\ldots.
\end{align}

In order to get the special form of $C$ in Kaczmarz-Tanabe's iteration, Kang introduced the concept of index set $I_d(n_1,n_2,v)$ in \cite[Definition 2.9]{Kang2023Thestandard}, i.e.,
\begin{defn}\label{definition.index.set} The index set $I_d(n_1,n_2,v)$ is defined as follows
\begin{align*}
  I_d(n_1,n_2,v)=\Big\{[I_d(1),\ldots,I_d(v)]\Big\},
\end{align*}
where $n_1,n_2,v$ are positive integers satisfied $|n_1-n_2|\ge v\ge 2$. $I_d(i)$ is integer between $n_1$ and $n_2$, and
$I_d(1)=n_1,I_d(v)=n_2$. For any $i<j$ there satisfy
\begin{align*}
\begin{array}{ll}
  I_d(i)<I_d(j), & n_1<n_2,\\
  I_d(i)>I_d(j), & n_1>n_2.
\end{array}
\end{align*}
\end{defn}

By Definition \ref{definition.index.set}, we can give the following expression of $a_i^T[Q_i({\bf u}_i)]^T\tilde{x}$ for any $1\le i\le m$ and $\tilde{x}\in N(A)^\bot$.
\begin{lem}\label{lemma.compatible.vector} Suppose that $A$ has no zero row, $Q_i$ is the sequential projection matrix of $A$ and $\tilde{x}\in N(A)^\bot$. For any $1\le i\le m, i+1\le j\le m-1$, denote
\begin{align}\label{C.element}
   d_{i,j}({\bf u}_i)=\sum_{v=2}^{j-i+1}(-1)^{v-1}\sum_{I_d(i,j,v)}\prod_{s=1}^{v-1}\mu_{I_d(s+1)}\prod_{s=1}^{v-1} h_{I_d(s),I_d(s+1)}.
\end{align}
Then,
\begin{align}\label{compatible.expression.1}
  a_i^T[Q_i({\bf u}_i)]^T\tilde{x}=(1,d_{i,i+1}({\bf u}_i),\ldots,d_{i,m}({\bf u}_i))(a_i^T,a_{i+1}^T,\ldots,a_m^T)^T\tilde{x}
\end{align}
holds. That is, the compatible vector of $a_i^T[Q_i({\bf u}_i)]^T$ on $A$ is $(0,\ldots, 0, 1,d_{i,i+1}({\bf u}_i),\ldots,d_{i,m}({\bf u}_i))$.
\end{lem}
\proof \quad When $1\le i\le m$ and $\tilde{x}\in N(A)^\bot$, we have
\begin{align*}
  a_i^T[Q_i({\bf u}_i)]^T\tilde{x}
  &=(1,-\mu_{i+1}h_{i,i+1})(a_i^T[Q_{i+1}({\bf u}_{i+1})]^T\tilde{x},a_{i+1}^T[Q_{i+1}({\bf u}_{i+1})]^T\tilde{x})^T\nonumber\\
  &=(1,-\mu_{i+1}h_{i,i+1},-\mu_{i+2}h_{i,i+2}+\mu_{i+1}\mu_{i+2}h_{i,i+1}h_{i+1,i+2})\nonumber\\
  & \quad \cdot (a_i^T[Q_{i+2}({\bf u}_{i+2})]^T\tilde{x},a_{i+1}^T[Q_{i+2}({\bf u}_{i+2})]^T\tilde{x},a_{i+2}^T[Q_{i+2}({\bf u}_{i+2})]^T\tilde{x})^T\nonumber\\
  &=\big(1,-\mu_{i+1}h_{i,i+1},\ldots,\sum_{v=2}^{m-i+1}(-1)^{m-1}\sum_{I_d(i,m,v)}\prod_{s=1}^{v-1}\mu_{I_d(s+1)}\prod_{s=1}^{v-1} h_{I_d(s),I_d(s+1)}\big)\nonumber\\
  &\quad \cdot (a_i^T[Q_m({\bf u}_m)]^T\tilde{x},a_{i+1}^T[Q_m({\bf u}_m)]^T\tilde{x},\ldots,a_m^T[Q_m({\bf u}_m)]^T\tilde{x})^T.
\end{align*}
By the agreement $Q_m({\bf u}_m)=I$,
\begin{align}\label{C.element.equivalent}
  a_i^T[Q_i({\bf u}_i)]^T\tilde{x}
  &=\big(1,-\mu_{i+1}h_{i,i+1},\ldots,\sum_{v=2}^{m-i+1}(-1)^{v-1}\sum_{I_d(i,m,v)}\prod_{s=1}^{v-1}\mu_{I_d(s+1)}\prod_{s=1}^{v-1} h_{I_d(s),I_d(s+1)}\big)\nonumber\\
  &\quad \cdot (a_i^T\tilde{x},a_{i+1}^T\tilde{x},\ldots,a_m^T\tilde{x})^T.
\end{align}
Thus \eqref{compatible.expression.1} holds by taking $d_{i,j}({\bf u}_i)$ according to \eqref{C.element}. Because \eqref{compatible.expression.1} holds for any $\tilde{x}\in \mathbb{R}^n$, we then get the compatible vector of $a_i^T[Q_i({\bf u}_i)]^T$ on $A$. \qed

\begin{thm}\label{thm.decomposition.detailed} Under the Lemma \ref{lemma.compatible.vector}, let $\Omega({\bf u})=(\omega_{i,j}({\bf u}_i))_{m\times m}$ satisfy
\begin{align}\label{definition.Omega.u}
  \omega_{i,j}({\bf u}_i)=\left \{
      \begin{array}{ll}
        d_{i,j}({\bf u}_i), & j>i,\\
        1,       & j=i,\\
        0      , & j<i.
      \end{array}
      \right.
\end{align}
Then,
\begin{align*}
  A_\mathcal{S}({\bf u})=\Omega({\bf u}) A.
\end{align*}
\end{thm}
\proof \quad For any $\tilde{x}\in N(A)^\bot$,
\begin{align}\label{As.x}
  A_\mathcal{S}({\bf u})\tilde{x}=(a_1^T[Q_1({\bf u}_1)]^T\tilde{x},a_2^T[Q_2({\bf u}_2)]^T\tilde{x},\ldots,a_m^T[Q_m({\bf u}_m)]^T\tilde{x})^T.
\end{align}
By Lemma \ref{lemma.compatible.vector} and \eqref{definition.Omega.u}, we then get
\begin{align}\label{A.s.equation}
  A_\mathcal{S}({\bf u})\tilde{x}=\Omega({\bf u}) A\tilde{x}.
\end{align}
Moreover, when $\tilde{x}\in N(A)$, \eqref{A.s.equation} obviously holds. Therefore, for any $\tilde{x}\in \mathbb{R}^n$, there holds
$A_\mathcal{S}({\bf u})\tilde{x}=\Omega({\bf u}) A\tilde{x}$, which means $A_\mathcal{S}({\bf u})=\Omega({\bf u}) A$. \qed

In the context, we especially call the iteration \eqref{normal.form.relaxation.Kaczmarz.Tanabe.iteration} \textbf{\emph{relaxation Kaczmarz-Tanabe's iteration}} when $C({\bf u})=\Omega({\bf u})$ is taken according to \eqref{definition.Omega.u}.

Let $E(j,i(-\mu_ih_{j,i}))$ be a matrix obtained by multiplying the $i$-th row of the identity matrix by $-\mu_ih_{j,i}$ and adding it to the $j$-th row, i.e., the diagonal elements of $E(j,i(-\mu_ih_{j,i}))$ are all $1$, the $(j,i)$- element is $-\mu_ih_{j,i}$, and all other elements are $0$.
Consequently, we have the following decomposition theorem of $\Omega({\bf u})$.

\begin{thm}\label{thm.Cmu.expression} For any $\mu_i\in \mathbb{R}$, let $\Omega({\bf u})$ be defined by \eqref{definition.Omega.u}, and define
\begin{align*}
   H_1(\mu_1)=I, H_i(\mu_i)=\prod\limits_{j=1}^{i-1} E(j,i(-\mu_i h_{j,i})),\quad 1<i\le m.
\end{align*}
Then,
\begin{align}\label{Cmu.decomposition}
   \Omega({\bf u})=H_1(\mu_1)H_2(\mu_2)\cdots H_m(\mu_m).
\end{align}
\end{thm}
\proof \quad For any $\tilde{x}\in N(A)^{\bot}$, from \eqref{C.element.equivalent}, we have
\begin{flalign}\label{formula.media.1}
  a_j^T[Q_j({\bf u}_j)]^T\tilde{x}
  &=\big(0,\ldots,0,1,-\mu_{j+1}h_{j,j+1},\ldots,\sum_{v=2}^{m-j+1}(-1)^{v-1}\sum_{I_d(j,m,v)}\prod_{s=1}^{v-1}\mu_{I_d(s+1)}\prod_{s=1}^{v-1} h_{I_d(s),I_d(s+1)}\big)\nonumber\\
  &\quad\cdot(a_1^T\tilde{x},a_2^T\tilde{x},\ldots,a_m^T\tilde{x})^T.
\end{flalign}
From \eqref{formula.media.1}, the coefficient of $a_i^T\tilde{x}$ is actually the $(j,i)$-element of $\Omega({\bf u})$ when $i>j$, i.e.,
\begin{align*}
   \omega_{j,i}({\bf u}_j)
   =-\mu_ih_{j,i}+\mu_i\mu_{i-1}h_{j,i-1}h_{i-1,i}+\ldots+(-1)^{i-j}\mu_i\mu_{i-1}\ldots\mu_{j+1}h_{j,j+1}h_{j+1,j+2}\cdot\ldots\cdot h_{i-1,i}.
\end{align*}
Denote $\widehat{H}({\bf u})=H_1(\mu_1)\cdots H_m(\mu_m)$, in order to show $\Omega({\bf u})=H_1(\mu_1)\cdots H_m(\mu_m)$, we only need to prove $\omega_{j,i}({\bf u}_j)=[\widehat{H}({\bf u})]_{j,i}$, i.e.,
\begin{align*}
  \omega_{j,i}({\bf u}_j)=e_j^T\widehat{H}({\bf u})e_i,
\end{align*}
where $e_j$ and $e_i$ are the $j$th and $i$th columns of the identity matrix in $\mathbb{R}^m$, respectively. Owing to $e_j^TH_k(\mu_k)=e_j^T$ when $j\ge k$ and $H_l(\mu_l)e_i=e_i$ when $i\neq l$, when $i>j$ it then follows from the notation of $\widehat{H}({\bf u})$ that
\begin{align*}
  e_j^T\widehat{H}({\bf u})e_i
  &=e_j^TH_{j+1}(\mu_{j+1})\cdots H_i(\mu_i)e_i\nonumber\\
  &=(e_j^TH_{j+1}(\mu_{j+1}))H_{j+2}(\mu_{j+2})\cdots H_i(\mu_{i})e_i\nonumber\\
  &=((0,\ldots,1,-\mu_{j+1}h_{j,j+1},0,\ldots,0)H_{j+2}(\mu_{j+2}))H_{j+3}(\mu_{j+3})\cdots H_i(\mu_i)e_i\nonumber\\
  &=((0,\ldots,1,-\mu_{j+1}h_{j,j+1},-\mu_{j+1}h_{j,j+1}+\mu_{j+2}\mu_{j+1}h_{j,j+1}h_{j+1,j+2},\ldots,0)H_{j+2}(\mu_{j+2}))\\
  & \quad \cdot H_{j+3}(\mu_{j+3})\cdots H_i(\mu_i)e_i\nonumber\\
  &=(0,\ldots,1,-\mu_{j+1}h_{j,j+1},\ldots,-\mu_ih_{j,i}+\mu_{i}\mu_{i-1}h_{j,i-1}h_{i-1,i}+\ldots\nonumber\\
  &\quad +(-1)^{i-j}\mu_{i}\mu_{i-1}\ldots\mu_{j+1}h_{j,j+1}\cdots h_{i-1,i},0,\ldots,0)e_i\nonumber\\
  &=-\mu_{i}h_{j,i}+\mu_{i}\mu_{i-1}h_{j,i-1}h_{i-1,i}+\ldots+(-1)^{i-j}\mu_{i}\mu_{i-1}\cdots\mu_{j+1}h_{j,j+1}h_{j+1,j+2}\cdots h_{i-1,i}.&
\end{align*}
This proves
\[\omega_{j,i}({\bf u}_i)=[\widehat{H}({\bf u})]_{j,i}\]
for any $1\le j\le n-1$ and $i>j$. Additionally, there also holds
\[\omega_{j,j}({\bf u}_j)=[\widehat{H}({\bf u})]_{j,j}=1\]
for any $1\le j\le n$. Consequently, the conclusion is proved. \qed

Theorem \ref{thm.Cmu.expression} actually provides a computational method of $C({\bf u})$ for the relaxation Kaczmarz-Tanabe's algorithm. Moreover, there also holds the following decomposition about $\Lambda C({\bf u})$.
\begin{lem}\label{lemma.Cu.decomposition} Let $C({\bf u})$ be defined like $\Omega({\bf u})$ in Theorem \ref{thm.decomposition.detailed}, then there exists decomposition of
\begin{align}\label{Cu.decomposition.formula}
  \Lambda C({\bf u})=T_u({\bf u})M,
\end{align}
where $T_u({\bf u})\in\mathbb{R}^{m\times m}$ is an upper triangular matrix with diagonal elements $[T_u({\bf u})](i,i)=\mu_i\|a_i\|_2^2$.
\end{lem}
\proof \quad The proof of the lemma is trivial and omitted here. \qed

\subsection{The convergence of the relaxation Kaczmarz-Tanabe method}

Next, we consider the convergence of the relaxation Kaczmarz-Tanabe method and give a sufficient condition of convergence for the relaxation Kaczmarz-Tanabe's iteration. Notice that each relaxation Kaczmarz-Tanabe's iteration is a composite of $m$ relaxation Kaczmarz's iterations, consequently, it is natural to transfer the convergence condition of the relaxation Kaczmarz method to the relaxation Kaczmarz-Tanabe method.

\begin{lem}\label{lemma.kernal.space.1} For a linear system $Ax=b$, let $Q({\bf u})$ be defined by \eqref{definition.Qmu}, and we suppose $0<\mu_i<2$ for each $1\le i\le m$. Then, for any $x\in \mathbb{R}^n$, there holds $ \|Q({\bf u})x\|_2\le \|x\|_2$, and the equality holds iif $x\in N(A)$.\end{lem}
\proof \quad When $x\in N(A)$, from $P_j(\mu_j)x=x, j=1,\ldots,m$, therefore there holds
\begin{align}\label{kernal.space.1}
  Q({\bf u})x=P_m(\mu_m)\cdots P_1(\mu_1)x=x.
\end{align}
If $x\notin N(A)$, from Lemma \ref{lemma.Pi.norm}, there exists at least an index $s_0\in\Theta\triangleq\{1,2,\ldots,m\}$ such that $\|P_{s_0}(\mu_{s_0})x\|_2<\|x\|_2$. We suppose that $s_0$ is the smallest index in $\Theta$ that $\|P_{s_0}(\mu_{s_0})x\|_2<\|x\|_2$, then there holds
\[\|Q({\bf u})x\|_2=\|P_m(\mu_m)\cdots P_1(\mu_1)x\|_2\le\|P_{s_0}(\mu_{s_0})\cdots P_1(\mu_1)x\|_2\le \|p_{s_0}(\mu_{s_0})x\|_2< \|x\|_2. \qedhere\]
These prove the conclusions. \qed

Lemma \ref{lemma.kernal.space.1} actually gives the condition to ensure the nonexpansion of $Q({\bf u})$ on $\mathbb{R}^n$.
Let $\{y_k\}$ be generated by \eqref{normal.form.relaxation.Kaczmarz.Tanabe.iteration} and $\bar{e}_k=y_k-x^\dagger-P_{N(A)}y_0$, then
\begin{align*}
  \bar{e}_k=\big(I-A^T[C({\bf u})]^T\Lambda MA\big)\bar{e}_{k-1},
\end{align*}
it follows from Lemma \ref{lem.Qmu} that
\begin{align}\label{error.relaxation.Kaczmarz.Tanabe}
  \bar{e}_k=Q({\bf u})\bar{e}_{k-1}.
\end{align}

By Lemma \ref{lemma.kernal.space.1}, thus we have the following convergence of the relaxation Kaczmarz-Tanabe method.
\begin{thm} \label{thm.convergence.relaxation.Kaczmarz.Tanabe} For a consistent linear system \eqref{linear.system}, supposing that $\{y_k,k\ge 0\}$ is generated by the relaxation Kaczmarz-Tanabe's iteration \eqref{normal.form.relaxation.Kaczmarz.Tanabe.iteration} with $0<\mu_i<2$ (where $1\le i\le m$). Then $\{y_k, k>0\}$ converges to $x^\dagger+P_{N(A)}y_0$.
\end{thm}
\proof \quad Let $Q_r({\bf u})=Q({\bf u})|_{N(A)^{\bot}}$, by Lemma \ref{lemma.kernal.space.1}, we have $\|Q_r({\bf u})\|_2<1$. Because the sequence generated by the relaxation Kaczmarz-Tanabe's iteration is a subsequence of the sequence generated by the relaxation Kaczmarz's iteration, there holds from Lemma \ref{lem.set.property.relaxation.Kaczmarz} that $\bar{e}_k\in N(A)^{\bot}$. Consequently, by \eqref{error.relaxation.Kaczmarz.Tanabe}, we have
\begin{align*}
  \|\bar{e}_k\|_2=\|Q({\bf u})\bar{e}_{k-1}\|_2=\|Q_r({\bf u})\bar{e}_{k-1}\|_2\le \|Q_r({\bf u})\|_2\|\bar{e}_{k-1}\|_2\le\cdots\le\|Q_r({\bf u})\|_2^k\|\bar{e}_0\|_2,
\end{align*}
which shows $\|\bar{e}_k\|\rightarrow 0$ as $k\rightarrow \infty$, i.e., $y_k\rightarrow x^\dagger+P_{N(A)}y_0$.\qed

Moreover, we have the following convergence rate of the relaxation Kaczmarz-Tanabe method for solving error-free linear system.
\begin{cor}\label{thm.convergence.rate.consistent} Under the conditions of Theorem \ref{thm.convergence.relaxation.Kaczmarz.Tanabe}, there also hold
\begin{align}\label{convergence.rate.relaxation.kt.1}
   \|\bar{e}_k\|_2\le\max_{0<\sigma_i({\bf u})<1}\sigma_i({\bf u})\|\bar{e}_{k-1}\|_2
\end{align}
and
\begin{align}\label{convergence.rate.relaxation.kt.2}
  \|\bar{e}_k\|_2\le\max_{0<\sigma_i({\bf u})<1}[\sigma_i({\bf u})]^k\|\bar{e}_0\|_2,
\end{align}
where $\sigma_i({\bf{u}})$ is the singular value of $Q({\bf u})$.
\end{cor}
\proof \quad From the proof of Theorem \ref{thm.convergence.relaxation.Kaczmarz.Tanabe}, we have
\begin{align*}
  \|\bar{e}_k\|_2\le\|Q_r({\bf u})\|_2\|\bar{e}_{k-1}\|_2.
\end{align*}
When $0<\mu_i<2$, $\|Q({\bf u})\|_2\le 1$ and $\|Q_r({\bf u})\|_2<1$. Then,
\begin{align*}
  \|Q_r({\bf u})\|_2\le \max_{0<\sigma_i({\bf u})<1}[\sigma_i({\bf u})].
\end{align*}
Consequently, \eqref{convergence.rate.relaxation.kt.1} and \eqref{convergence.rate.relaxation.kt.2} hold. \qed

\section{The related algorithms}
\label{section.algorithm}

The implementation of the relaxation Kaczmarz-Tanabe method are divided into two processes. The first process is to compute $C({\bf u})$ for given ${\bf u}\in \mathbb{R}^m$, and the second process is to execute the Kaczmarz-Tanabe's iteration repeatedly. In Algorithm \ref{algorithm.C.flow}, we list the flowchart to compute $C({\bf u})$. After computing $C({\bf u})$, the implementation of the relaxation Kaczmarz-Tanabe's iteration is trivial, and we list the process in Algorithm \ref{algorithm.relaxation.KT.iteration}.

\begin{breakablealgorithm}\label{algorithm.C.flow}
\caption{The calculation of matrix $C({\bf u})$}
\begin{algorithmic}[1]
\State \textbf{Input}
\State \indent $A=(a_1,a_2,\ldots,a_m)^T$
\State \indent ${\bf u}=(\mu_1,\mu_2,\ldots,\mu_m)$  \Comment{${\bf u}$ relaxation parameter, $\mu_i=1$ is for decomposition without relaxation}
\State \indent $C=I_m$ \Comment{$I_m$ is an identity matrix with order $m$}
\State \indent $k \gets m$
\While{$k>1$}
  \State $i \gets k$
  \While{$i>1$}
     \State $j \gets m$
     \While{$j>k-1$}
     \If {$a_k^Ta_k=0$}
        \State $C(i-1,j)=C(i-1,j)$
     \Else
        \State $C(i-1,j)=C(i-1,j)+\mu_k(-a_{i-1}^Ta_k/(a_k^Ta_k)C(k,j)$
     \EndIf
        \State $j \gets j-1$
     \EndWhile
  \State $i \gets i-1$
  \EndWhile
\State $k \gets k-1$
\EndWhile
\State \textbf{Output} $C({\bf u})$
\end{algorithmic}
\end{breakablealgorithm}

\begin{breakablealgorithm}
\label{algorithm.relaxation.KT.iteration}
\caption{The relaxation Kaczmarz-Tanabe method solver}
\begin{algorithmic}[1]
\State \textbf{Input}
\State \indent $A,M,\Lambda,b,{\bf u}$
\State \indent $y_0=x_0$ \Comment{Initial value}
\State \indent $K_{\max}$ \Comment{Maximal iteration number}

\State \textbf{Compute} $C({\bf u})$ by Algorithm \ref{algorithm.C.flow}
\State $k=1$
\While{$k<K_{\max}$}
  \State $y_{k+1}=y_k+A^T[C({\bf u})]^T\Lambda M(b-Ay_k)$
  \State $k\gets k+1$
\EndWhile
\State \textbf{Output} $y_k$. \Comment{$y_k$ is the acquired solution}
\end{algorithmic}
\end{breakablealgorithm}

\section{Conclusion}
\label{section.conclusion}

In this paper, we consider the relaxation Kaczmarz-Tanabe method by introducing the relaxation parameter idea into Kaczmarz-Tanabe method. The development of this work is based on the research to the Kaczmarz-Tanabe method. We first make some exploration on the selections of the relaxation parameter for the relaxation Kaczmarz-Tanabe method, then analyze the related theories of relaxation Kaczmarz-Tanabe's iteration such as the the decompositions of $A_{\mathcal{S}}({\bf u})$ and $C({\bf u})$.

In numerical tests, we experiment several classic examples to verify these theoretical results of the selections of relaxation parameters for the relaxation Kaczmarz-Tanabe method. Although we prove that $\mu_i\in (0,2)$ is the sufficient condition for the convergence of the relaxation Kaczmarz-Tanabe method, the operation of adjusting parameters one by one is tedious and unnecessary. Consequently, we carry out many experiments with constant and random parameters to test the effect of parameter selection, moreover the experimental results of random parameters setting show the effectiveness of the specific case. However the selection of optimal parameters still needs further research, and in the current situation $\mu_i\equiv 1$ may still be a compromise for the relaxation Kaczmarz-Tanabe method.

\section*{Reference}

\bibliographystyle{unsrt}
\bibliography{Reference}

\end{document}